\documentclass{article}

\usepackage{PRIMEarxiv}
\usepackage{authblk}
\usepackage{graphicx,xcolor,enumerate,amssymb,amsmath,verbatim}
\usepackage{amssymb}
\usepackage{amsbsy}
\usepackage{float}
\usepackage{epstopdf}
\usepackage[round, sort]{natbib}
\usepackage[utf8]{inputenc} % allow utf-8 input
\usepackage[T1]{fontenc}    % use 8-bit T1 fonts
\usepackage{hyperref}       % hyperlinks
\usepackage{url}            % simple URL typesetting
\usepackage{booktabs}       % professional-quality tables
\usepackage{amsfonts}       % blackboard math symbols
\usepackage{nicefrac}       % compact symbols for 1/2, etc.
\usepackage{microtype}      % microtypography
\usepackage{lipsum}
\usepackage{fancyhdr}       % header
\usepackage{graphicx}       % graphics
\graphicspath{{media/}}     % organize your images and other figures under media/ folder

\makeatletter
\@addtoreset{equation}{section}
\renewcommand\thetable{\@arabic\c@table}
\renewcommand\thefigure{\@arabic\c@figure}
\long\def\@makecaption#1#2{%
  \vskip\abovecaptionskip
  \begin{center}%
  \sbox\@tempboxa{#1: #2}%
  \ifdim \wd\@tempboxa >\hsize
    #1: #2\par
  \else
    \global \@minipagefalse
    \hb@xt@\hsize{\hfil\box\@tempboxa\hfil}%
  \fi
  \end{center}%
  \vskip\belowcaptionskip}
\def\N{{\rm I\kern-.15em N}}
\def\R{{\rm I\kern-.2em R}}
\def\Z{{\rm Z\kern-.26em Z}}
\makeatother

\newtheorem{thm}{Theorem}[section]

\newcommand{\be}{\begin{eqnarray}}
\newcommand{\ee}{\end{eqnarray}}
\newcommand{\bq}{\begin{eqnarray*}}
\newcommand{\eq}{\end{eqnarray*}}

\newcommand{\stk}{\stackrel{\mathbb{\scriptsize P}}{\longrightarrow}}

\newcommand{\verk}{\stackrel{{\cal D}}{\longrightarrow}}

\title{Bahadur efficiency for certain goodness--of--fit tests based on the empirical characteristic function
}

\author[a,b]{Simos G. Meintanis\thanks{simosmei@econ.uoa.gr}}
\affil[a]{\small Department of Economics, National and Kapodistrian University
of Athens, Athens, Greece}
\author[c]{Bojana Milo\v{s}evi\'c\thanks{bojana@matf.bg.ac.rs}}
\author[c]{Marko Obradovi\' c\thanks{marcone@matf.bg.ac.rs}}
\affil[b]{\small Pure and Applied Analytics\\ North--West University\\ Potchefstroom, South Africa}
\affil[c]{\small Faculty of Mathematics, University of Belgrade, Studenski trg 16, Belgrade, Serbia}
\date{}

\begin{document}
\maketitle

\begin{abstract}
We study the Bahadur efficiency of several weighted L2--type goodness--of--fit tests  based on the empirical characteristic function. The methods considered are for normality and exponentiality testing, and for testing goodness--of--fit to the logistic distribution. Our results are helpful in deciding which specific test a potential practitioner should apply.  For the celebrated BHEP and energy tests for normality we obtain novel efficiency results, with some of them in the multivariate case, while in the case of the logistic distribution this is the first time that efficiencies are computed for any composite goodness--of--fit test.
\end{abstract}

% keywords can be removed
\keywords{Goodness--of--fit test\and Bahadur efficiency\and Empirical characteristic function \and  Normality test \and Exponentiality test}
\section{Introduction}\label{sec_1}

% \section{General formulation} \label{sec_1}
Let $X_1,...,X_n$, denote independent copies of an arbitrary random variable $X\in\mathbb R^p$, and consider the problem of testing the composite goodness--of--fit (GOF) null hypothesis, \begin{equation}\label{null} {\cal{H}}_{0}: \mbox{The law of}\ X\in   {\cal{F}}_{\vartheta}, \ {\rm{for \ some}} 
\ \vartheta, \end{equation} where ${\cal{F}}_{\vartheta}:=\{F_{\vartheta}, \vartheta \in \Theta\}$ denotes a class of distributions
indexed by the parameter  $\vartheta\in\Theta\subseteq \mathbb R^q, \ q\geq 1$.

For certain popular distributions, such as the normal and the exponential distribution, there exist many GOF tests, while for others such as the logistic distribution the range of methods available is not so extended. In either case though a potential practitioner would like to have some quality measure on the basis of which one could choose amongst the existing GOF methods. In this connection, one of the most popular methods for test comparison is the so--called Bahadur efficiency that allows to compare the efficiency of any given GOF test vis--\'a--vis its optimal counterpart, which is the likelihood ratio (LR) test for the distribution under test against a specified alternative.

%As it is typical in GOF testing we compare a nonparametric estimator of the DF or CF 
%with a corresponding parametric estimator of the same quantity reflecting the null hypothesis. %Specifically, let
%\begin{equation}\label{EDF} F_n(x)=\frac{\{\# X_j's\leq x\}}{n}, \end{equation} 
%and
%\begin{equation} \label{ECF}  \varphi_n(t)=\frac{1}{n} \sum_{j=1}^n e^{i t X_j},
%\end{equation}
%be the empirical DF and the empirical CF, respectively. Also denote by ${\widehat{\boldsymbol{\vartheta}}}_n=(\widehat \alpha_n,\widehat \beta_n,\widehat c_n,\widehat %\delta_n)'$ the ML estimator of the parameter  
%${\boldsymbol{\vartheta}}=(\alpha,\beta,c,\delta)'$ and write $F_{{\widehat{\boldsymbol{\vartheta}}}_n}(\cdot)$ and $\varphi_{{\widehat{\boldsymbol{\vartheta}}}_n}(\cdot)$
%for the corresponding DF and CF of a SP law with the parameter ${\boldsymbol{\vartheta}}$ replaced by this estimator.  Then empirical DF GOF statistics may generally be expressed as \begin{equation}\label{DFTEST} D_n=\Delta(F_n,F_{{\widehat{\boldsymbol{\vartheta}}}_n}), \end{equation}  where  $\Delta(\cdot,\cdot)$, denotes some distance function, and likewise empirical CF GOF statistics may be written as   
%\begin{equation}\label{CFTEST} C_n=\Delta(\varphi_n,\varphi_{{\widehat{\boldsymbol{\vartheta}}}_n}), \end{equation} 
%and we reject the null hypothesis (\ref{null}) for large values of these test statistics. 

%\subsection{Goodness--of--fit tests based on the CF}
%\vspace*{-.3cm}
In this work we consider GOF tests that are based on the weighted L2--type test statistic
\begin{equation} \label{ts}
T_{n,w}=n\: \int_{\mathbb R^p} \left|\varphi_n(t)-\varphi_{\widehat\vartheta_n}(t)\right|^2 w(t) {\rm{d}}t,
\end{equation}
where $\varphi_\vartheta(\cdot)$ denotes the characteristic function corresponding to $F_\vartheta$,    
\begin{equation} \label{ECF}  \varphi_n(t)=\frac{1}{n} \sum_{j=1}^n e^{{\rm{i}} t^\top X_j},
\end{equation}
is the empirical characteristic function (ECF) and $\widehat\vartheta_n:=\widehat\vartheta_n(X_1,...,X_n)$ is an
estimator of $\vartheta$ obtained on the basis of $(X_j, \ j=1,...,n)$. 

The test statistic $T_{n,w}$, besides the family ${\cal {F}}_\vartheta$ being tested and the corresponding 
estimator $\widehat\vartheta_n$ employed, also depends on the the weight function $w\geq 0$ figuring in \eqref{ts}.
For certain choices of $w(\cdot)$, we obtain here efficiency results for the BHEP GOF test for normality of \cite{epps} and \cite{henze1997new} as well as
for the generalized energy (GE) test for normality first put forward by \cite{szekely2005new}. In doing so we provide extension of the efficiency results of \cite{ebner2021bahadur} and \cite{tenreiro2009choice} for the BHEP test as well as extension of the results of  \cite{mori2021energy} from the original energy statistic to its generalized counterpart suggested by \cite{szekely2013energy}. Furthermore we obtain for the first time analogous results for the GOF test for the logistic distribution of Meintanis (2004) and for the exponentiality test suggested  by Henze and Meintanis (2005), including efficiency comparisons with alternative tests. As already mentioned these efficiency results will facilitate the choice of the specific weight function $w(\cdot)$ that figures in the test statistics defined by \eqref{ts}, and thus provide some guidance on which method to apply amongst the many tests available for the normal, the exponential distribution and the logistic distribution.     

The rest of this work unfolds as follows. In Section \ref{sec2} we review the basic theory of Bahadur efficiency and the theory of the limit null distribution of the test statistics, while in Section \ref{sec3} the necessary efficiency computations are discussed. In Section \ref{sec4} we provide efficiency results and discussion for the BHEP and energy statistics, as well as for the two aforementioned ECF--based tests, one for the logistic and the other for the exponential distribution. Analogous efficiency results for the normality tests in the bivariate case are reported in Section \ref{sec5}, and we conclude in Section \ref{sec6} with discussion.

\section{Limit null distribution and Bahadur slopes} \label{sec2}

%\section{Limit null distribution and Bahadur slope}
The limit in distribution of the L2--type test statistics figuring in 
\eqref{ts} is given by    
\be \label{Tw}
T_{n,w} \verk \sum_{j=1}^\infty \lambda_j N^2_j, 
\ee
where $N_j, \ j\geq 1$, are independent copies of a standard normally distributed random variable, 
and $\lambda_1\geq \lambda_2\geq \ldots$,  are the eigenvalues of the integral equation 
\be \label{eigen}
\lambda f(s)=\int K(s,t) f(t) w(t) {\rm{d}}t;  
\ee where $K(s,t)$ is a covariance kernel associated, besides the weight function $w(\cdot)$,  with the given test statistic and the estimator employed in estimating the distributional parameter $\vartheta$; see  for instance \cite{meintanis2007bootstrap}.  

In this section, we briefly review the essence of the Bahadur theory, and associate that theory with the limit law of $T_{n,w}$ figuring in \eqref{Tw}--\eqref{eigen}; for more details we refer to \cite{bahadur1971some} and \cite{nikitinKnjiga}, and to the more recent article by  \cite{grane2}. To this end let $\mathcal{G}_\theta=\{G_\theta,\;\theta>0\}$ be a family of alternative distribution functions, 
such that $G_0$ is the null family of distributions for some typical value $\vartheta\in \Theta$. Assume that the regularity conditions for 
V--statistics with weakly degenerate kernels in \cite{nikitinMetron},  are satisfied. 

Also recall that LR tests are optimal tests in the Bahadur sense.
Hence, for close alternatives from $\mathcal{G}_\theta$,  the absolute local approximate Bahadur efficiency for any 
sequence  of test statistics, say  $\{T_n\}$, is defined as the ratio of the Bahadur approximate slope of the 
considered test statistic to the corresponding slope of the LR test, i.e.
\begin{equation}\label{eff}
 {\rm eff}(T_n)=\lim\limits_{\theta\to 0}\frac{c_T(\theta)}{2K(\theta)},
\end{equation}
where \be \label {cT}
c_{T}(\theta)=\frac{b_T(\theta)}{\lambda_1}=\frac{b''_T(0)}{2\lambda_1}\theta^2+{\rm{o}}(\theta^2),
\ee
(the last equation obtained by a Taylor expansion), with $\lambda_1$ being the largest eigenvalue figuring in \eqref{Tw} and 
\be \label{beta}
\frac{T_n}{n} \stk b_T(\theta),
\ee
where $K(\theta)$ is equal to minimal Kullback--Leibler (KL) 
distance from the given alternative to the class of distributions within the null hypothesis.

%In this connection, it is well known that the celebrated measure of efficiency called local Bahadur slope 
%critically depends on the largest eigenvalue $\lambda_{\max}$.
%In the univariate case, \cite{tenreiro2009choice}  has obtained eigenvalue--approximations for the 
%BHEP test of normality, while  \cite{mori2021energy} studied the same problem for the energy test first 
%put forward by \cite{szekely2005new}. The eigenvalue approximation program for ECF--based tests in the case 
%of the more general stable family of distributions has been carried out by \cite{matsui2008goodness}.       

Before  closing  this section we wish to point out that the BHEP and GE tests are well known to be location/scale invariant, i.e. for both  tests it holds,  \begin{equation} \label{affine}
T_{n,w}(a+bX_1,...,a+bX_n)=T_{n,w}(X_1,...,X_n), \ \mbox{for \ each}\  a\in \mathbb R, b>0,
\end{equation}
and therefore their respective limit null distributions are independent of the true values of the mean 
and variance of the underlying Gaussian law. Note that this property extends to the multivariate version of these tests. In this connection, the other tests studied herein are also invariant, whereby invariance is understood within the context of the specific family ${\cal{F}}_\vartheta$ being tested. For instance if this family is confined to the positive real axis (such is the exponential family of distributions), invariance  is understood only with respect to scale, i.e.   \eqref{affine} holds for $a=0$, and $b>0$. %We refer to \cite{henze2002invariant} for a detailed discussion of affine invariannt normality tests in arbitrary dimension, and to  \cite{meintanis2007bootstrap} for the invariance properties of
%ECF--based test statistics.    

\section{KL distance and eigenvalue approximation} \label{sec3} 
It can be seen from \eqref{eff}--\eqref{beta}   that the basic components of an approximate Bahadur slope are the computation of the KL distance and of the largest eigenvalue. In this section we consider these components with emphasis on location--scale families of distributions.    
\subsection{KL distance for location-scale families} 
Assume that the family $\mathcal{F}_\vartheta$ is a location--scale family, i.e. a family  generated by the density function $f_\vartheta$ corresponding to $F_\vartheta$, where $\vartheta$ contains the location and  scale parameters, $\mu$ and $\sigma$, respectively. We set  $f_0(\cdot)$ and $F_0(\cdot)$ for the density and distribution function, respectively, of the standardized variable $Z:=\sigma^{-1}(X-\mu)$. Recall also that for the family $\mathcal{G}_\theta$ with density $g_\theta(\cdot), \; \theta>0$, we have that $g_0\equiv f_0$,  and assume that certain regularity conditions are satisfied. Then the next theorem gives the local behavior of the KL distance $K(\theta;\mu,\sigma)$ between $g_\theta\in\mathcal{G}_\theta$ and $f_0 \in {\cal{F}}_\vartheta$,  
\begin{align}\label{KL}
    K(\theta;\mu,\sigma)&=\int\log(g_\theta(x))g_\theta(x){\rm{d}}x-\int\log(f_0(x))g_\theta(x){\rm{d}}x. 
\end{align}
The proof of Theorem \ref{TKL} is postponed to the Appendix.

\begin{thm}\label{TKL}
%\begin{thm}{Theorem}[section]

    The KL distance from the alternative $g_\theta(\cdot)$ in $\mathcal{G}_\theta$ to the closest null distribution
    \begin{align*}K(\theta)=\inf_{\mu,\sigma}K(\theta;\mu,\sigma)
    \end{align*}
    admits the representation
    \begin{align} \label{KLexpr} 2K(\theta)&=\Bigg(\int\frac{h^2(x)}{f_0(x)}{\rm{d}}x-(\sigma'(0))^2+\int\frac{1}{f_0(x)}(f_0'(x))^2(\mu'(0)+x\sigma'(0))^2{\rm{d}}x+\\ \nonumber &2\int\frac{1}{f_0(x)}f_0'(x)(\mu'(0)+x\sigma'(0))h(x){\rm{d}}x\Bigg) \theta^2 + o(\theta^2), \;\;\theta\to 0,
    \end{align}
    where $(\mu(\theta),\sigma(\theta))={\rm arginf} \: K(\theta;\mu,\sigma)$ and $h(x)=\frac{\partial}{\partial \theta}g_{\theta}(x)|_{\theta=0}$. %(All derivatives are with respect to $\theta$). 
   \end{thm}

   In the case of some specific distribution the expression above can be simplified. Since the maximum likelihood (ML)
   estimators minimize the KL distance, the functions $\mu(\theta)$ and $\sigma(\theta)$ are equal to population counterparts, 
   i.e. probability limits of the ML estimators of the location and scale parameter, respectively.
   
   For example, in the case that ${\cal{F}}_\vartheta$ is the normal location--scale family, we have 
   \begin{align*}\mu(\theta)&=\int_{-\infty}^\infty xg_\theta(x){\rm{d}}x; \ \ 
    \sigma(\theta)= \bigg(\int_{-\infty}^\infty x^2g_\theta(x){\rm{d}}x - \mu^2(\theta)\bigg)^{1/2}
   \end{align*}
   and thus the KL distance reduces to (see also \cite{milovsevic2021bahadur})
   \begin{align*}
    2K(\theta)&=\Bigg(\int_{-\infty}^{\infty}\sqrt{2\pi}h^2(x)e^{\frac{x^2}2}{\rm{d}}x-
    \Big(\int_{-\infty}^{\infty}xh(x){\rm{d}}x\Big)^2\\&-
    \frac12\Big(\int_{-\infty}^{\infty}x^2h(x){\rm{d}}x\Big)^2\Bigg)\theta^2 + o(\theta^2).
           \end{align*}

For the exponential scale family, we have
   \begin{align*}\mu(\theta)=0; \ \     \sigma(\theta)=\int_0^\infty xg_\theta(x){\rm{d}}x, 
   \end{align*}
%   and thus the KL distance reduces to (see also \cite{milovsevic2021bahadur})
and thus the KL distance reduces to (see also \cite{nikitin1996bahadur})
\begin{align*}
 2K(\theta)=\bigg(\int_0^\infty h^2(x)e^x{\rm{d}}x -\Big(\int_0^\infty xh(x){\rm{d}}x\Big)^2\bigg) \theta^2+o(\theta^2).
\end{align*}

However, in the case that ${\cal{F}}_\vartheta$ is the logistic location--scale, the parameters $(\mu(\theta),\sigma(\theta))$ that minimize 
the KL distance figuring in \eqref{KL} admit  no closed--form expression. Nevertheless, by applying the implicit function theorem we obtain 
   \begin{align*}
   \mu'(0)=6 \int_{-\infty}^{\infty}\frac{h(x)}{1+e^{-x}}{\rm{d}}x;  \ \ 
\sigma'(0)=\frac{9}{\pi^2+3}\int_{-\infty}^{\infty}\frac{(1-e^{-x})xh(x)}{1+e^{-x}} {\rm{d}}x,
\end{align*}
and thus by plugging these expressions in \eqref{KLexpr} we can obtain the corresponding KL distance.

\subsection{Eigenvalue approximation}\label{sect:eigen}
In this section we briefly review a method for eigenvalue approximation proposed in \cite{bozin2020new}. Recall that these eigenvalues are solutions of an integral equation involving a specific operator. At first we replace the original operator by a symmetric operator that has the same eigenvalues, and then (i) consider a truncated version of the symmetrized operator, and (ii) consider a discretized version of the truncated operator. In these two steps the truncated operator is chosen so that as the amount of truncation diminishes this operator converges to the symmetrized operator, and likewise, the discretized operator approaches the truncated operator as the grid of discretization becomes more fine and extended.     

%The idea behind it is to express the largest eigenvalue as the limit of a sequence of largest eigenvalues of discrete linear operators that converge in norm to the integral operator in question; and then approximate it with a member of the sequence with arbitrarily large index. From the perturbation theory we know that this works if the integral operator (as well as discrete linear operators) is symmetric and self-adjoint. 

Specifically  we replace the original operator  
\begin{align*}
    Af(s)=\int_{-\infty}^{\infty}K(s,t)f(t)w(t){\rm{d}}t
\end{align*}
with the operator 
  \begin{equation}\label{operatorS}
  \begin{aligned}\overline{A}f(s)&=\int_{-\infty}^{\infty}K(s,t)f(t)\sqrt{w(t)w(s)}{\rm{d}}t,\end{aligned}
  \end{equation}
  that has the same spectrum as $A$, but is symmetric. Then in the first step we define the truncated operator $\overline{A_B}$ acting on the set of real functions with support $[-B,B], \ B>0$, defined  by   \begin{align*}
      \overline{A_B}f(s)=\int\limits_{-\infty}^{\infty}K(s,t)f(t)\sqrt{w(t)w(s)}\:{\rm 1}(|t|\leq B){\rm{d}}t,
  \end{align*}
which clearly and for sufficiently large $B$, %e.g. $B> -\log\varepsilon$ 
  is close to $\overline{A}$.

In the second step of approximation we employ a sequence of symmetric linear operators $M^{(m)}$, which converges in norm to $\overline{A_B}$, as $m\to \infty$. This discretized sequence  can be defined by $(2m+1)\times(2m+1)$ matrices $M^{(m)}=||m_{i,j}^{(m)}||,\; -m\leq i\leq m, -m\leq j\leq m$, with elements
\begin{equation}\label{MatAppr}
m_{i,j}^{(m)}=\frac{2B}{(2m+1)}K\left(\frac{iB}{m},\frac{jB}{m}\right)\sqrt{w\left(\frac{iB}{m}\right) w\left(\frac{jB}{m}\right)}.
\end{equation}

Using the perturbation theory--see \cite[Theorem 4.10, page 291]{kato2013perturbation}--we have that the spectra of these two operators are at a distance that
tends to zero. Hence within the degree of approximation, the sequence $\lambda_1^{(m)}$ of the largest eigenvalues of $M^{(m)}$ will converge to the largest eigenvalue $\lambda_1(B)$ of $\overline{A_B}$,  which in turn approaches $\lambda_1$. Consequently, the eigenvalues $\lambda_1^{(m)}$ and $\lambda_1$ will coincide up to any desired accuracy, provided that the pair of approximation parameters $(m,B)$ is large enough.

\section{Tests for the normal and the logistic distribution} \label{sec4}

%\subsection{Nonparametric alternatives}

When testing for GOF to a symmetric distribution such as the normal or the logistic, it is customary to use general purpose nonparametric 
alternatives parametrized by $\theta$, a parameter that usually controls skewness. 
A discussion on the construction and applications of these alternatives is available in \cite{jones2015families} and \cite{ley2015flexible}.  
Specifically we consider the following alternatives:

\begin{itemize}
\item Lehmann alternatives 
\begin{equation*}
        g^{(1)}_\theta(x)=(1+\theta)F_0^{\theta}(x)f_0(x)
    \end{equation*}
    \item first Ley-Paindaveine alternatives 
    \begin{equation*}
        g^{(2)}_\theta(x)=f_0(x) e^{-\theta(1-F_0(x))}(1+\theta F_0(x))
    \end{equation*}
    \item second Ley-Paindaveine alternatives
    \begin{equation*} 
        g^{(3)}_\theta(x)=f_0(x)(1-\theta\pi\cos(\pi F_0(x))
    \end{equation*}
    \item contamination alternatives 
    \begin{equation*}
        g^{(4)}_\theta(x;\mu,\sigma^2)=(1-\theta)f_0(x)+\theta f(x;\mu,\sigma^2) \label{contamination}
    \end{equation*}
    %\item Skew normal (Azzalini)
    %\begin{align}
        %g_5(x;\theta)=2\Phi(\theta x)\varphi(x)
    %\end{align}
\end{itemize}

Note that in the above alternatives, $f_0(\cdot)$ and $F_0(\cdot)$ are the null density and distribution function, respectively, so that for $\theta=0$ each of these alternatives reduces to the null distribution under test.

%\section{Tests based on the characteristic function}

\subsection{Tests for normality}

\subsubsection{The GE test} 
The GE test may be  formulated as 
\begin{equation} \label{energy}
T_{n,w}=n\: \int_{-\infty}^\infty \left|\phi_n(t)-e^{-\frac{t^2}{2}}\right|^2 w(t) {\rm{d}}t,
\end{equation}
with the ECF $\phi_n(t)$ obtained as in \eqref{ECF} with $X_j$ replaced by 
\be \label{stand} Z_j=\frac{X_j-\bar X_n}{S_n},\ j=1,...,n,\ee
where  $\bar X_n$ denotes the sample mean and $S^2_n$ the sample variance of $(X_j, \ j=1,...,n)$. As weight function \cite{mori2021energy} 
use $w(t)=|t|^{-2}$. Here we consider the GE test of \cite{szekely2013energy} whereby $w(t)=|t|^{-1-\gamma}, \ 0<\gamma<2$, 
is adopted as weight function, thus rendering the convenient form \be \label{energy1}
T_{n,w}=2 \sum_{j=1}^n \mathbb E|Z_j-X_1|^\gamma- n\mathbb E|X_1-X_2|^\gamma -\frac{1}{n} \sum_{j,k=1}^n |Z_j-Z_k|^\gamma.
\ee
\footnote {The original weight function also includes the 
constant $2\sqrt{\pi} \Gamma(1-(\gamma/2))/(\gamma 2^\gamma \Gamma((1+\gamma)/2))$} Note that  the 
expectations in \eqref{energy1} are taken with respect to the standard Gaussian distribution, and consequently the GE test statistic may be explicitly expressed by using the equations   
 \[\mathbb E|X_1-X_2|^\gamma=\frac{2^\gamma}{\sqrt{\pi}} \Gamma\left(\frac{1+\gamma}{2}\right),\]
%  \end{document}
and
 \begin{eqnarray*} \mathbb E|x-X_1|^\gamma&=&\frac{2^{\frac{\gamma}{2}}}{\sqrt{\pi}}
 e^{-\frac{x^2}{2}} \Gamma\left(\frac{1+\gamma}{2}\right) \: _{1}F_{1}\left(\frac{1}{2}+\frac{\gamma}{2},\frac{1}{2},\frac{x^2}{2}\right)\\ &=&\frac{2^{\frac{\gamma}{2}}}{\sqrt{\pi}}
 \Gamma\left(\frac{1+\gamma}{2}\right) \: _{1}F_{1}\left(-\frac{\gamma}{2},\frac{1}{2},-\frac{x^2}{2}\right),
\end{eqnarray*}
where $_{1}F_{1}(a,b,c)$ stands for the Kummer confluent hypergeometric function; see \cite{gradshteyn1994tables}.
For the covariance kernel $K(s,t)$   corresponding 
to the GE statistic  $T_{n,w}$ we refer to  \cite{mori2021energy}, 
while calculation of the eigenvalues  $\lambda_j, \ j\geq 1$, figuring in \eqref{Tw}--\eqref{eigen} 
is carried out using method described in \S\ref{sect:eigen}.

Letting $(\mu(\theta),\sigma^2(\theta))$ be the probability limit of the estimator $(\bar{X_n},\bar{S}^2_n)$, i.e.,
  \begin{align*}
      \mu(\theta)=\int_{-\infty}^{\infty}xg_\theta(x)dx; \ \ 
      \sigma^2(\theta)=\int_{-\infty}^{\infty}(x-\mu(\theta))^2g_\theta(x)dx,
  \end{align*}
we obtain the probability limit figuring in  \eqref{beta} for the GE test statistic as
\begin{align} \label{bTenergy} \frac{T_{n,w}}{n}\stk b(\theta)&=\frac{2^{1+\frac{\gamma}{2}}}{\sqrt{\pi}}\Gamma\left(\frac{1+\gamma}{2}\right)\mathbb E_{\theta}\left[
  \: _{1}F_{1}\left(-\frac{\gamma}{2},\frac{1}{2},-\frac{(X_1-\mu(\theta))^2}{2\sigma^2(\theta)}\right)\right]\\ \nonumber &-\frac{2^\gamma}{\sqrt{\pi}} \Gamma\left(\frac{1+\gamma}{2}\right)-\frac{1}{\sigma(\theta)^\gamma}\mathbb E_{\theta}\left|X_1-X_2\right|^{\gamma}.
  \end{align}
  %$b_1(\theta)=\int
  %\: _{1}F_{1}\left(-\frac{\gamma}{2},\frac{1}{2},-\frac{(X_1-\mu(\theta))^2}{2\sigma^2(\theta)}\right)g(x;\theta)dx$
  
  %\begin{align*}b'_1(0)&=\int -\frac{\gamma}{2} x \left(2  \mu'(0)+x \frac{d\sigma^2(0)}{d\theta}\right) \, _1F_1\left(1-\frac{\gamma}{2};\frac{3}{2};-\frac{x^2}{2}\right)g(x;0)dx\\&+
  %\int
  %\: _{1}F_{1}\left(-\frac{\gamma}{2},\frac{1}{2},-\frac{x^2}{2}\right)g'_{\theta}(x;0)dx
   % \end{align*}
 %$b_2(\theta)=\frac{1}{(\%sigma(\theta))^\gamma}\mathbb E_{\theta}\left|X_1-X_2\right|^{\gamma}$ 
 % \begin{align*}
  %    b_2'(0)= 2\int \frac{2^{\frac{\gamma}{2}}}{\sqrt{\pi}}
 %\Gamma\left(\frac{1+\gamma}{2}\right) \: _{1}F_{1}\left(-\frac{\gamma}{2},\frac{1}{2},-\frac{x_2^2}{2}\right)g'(x_2,0)dx_2-\frac{\gamma}{2}\frac{d\sigma^2(0)}{d\theta}\cdot \frac{2^\gamma}{\sqrt{\pi}} \Gamma\left(\frac{1+\gamma}{2}\right)
 % \end{align*}
  
 % \begin{align*}
  %    b'(0)&= \frac{2^{\frac{\gamma}{2}}}{\sqrt{\pi}}\Gamma\left(\frac{1+\gamma}{2}\right)  \bigg[ \int -\gamma x \left(2  \mu'(0)+x \frac{d\sigma^2(0)}{d\theta}\right) \, _1F_1\left(1-\frac{\gamma}{2};\frac{3}{2};-\frac{x^2}{2}\right)g(x;0)dx\\&\color{red}+\frac{\gamma 2^{\frac{\gamma}{2}}}{2}\frac{d\sigma^2(0)}{d\theta}\bigg]=0
  %\end{align*}

In order to calculate the quantity $b''(0)$ figuring in \eqref{cT} we use numerical integration in Wolfram Mathematica facilitated by
 differentiation under the integral sign, and thereby obtain the expressions 
  \begin{align*}
    \mu'(0)&=\int_{-\infty}^{\infty} x h(x)dx; \ \     \mu''(0)=\int_{-\infty}^{\infty} x u(x)dx\\
    \sigma'(0)&=\frac{1}{2}\int_{-\infty}^{\infty} x^2 h(x)dx\\
    \sigma''(0)&=\frac12\int_{-\infty}^{\infty} x^2 u(x)dx-\Big(\int_{-\infty}^{\infty} x h(x)dx\Big)^2 -\frac14\Big(\int_{-\infty}^{\infty} x^2 h(x)dx\Big)^2, 
\end{align*}
(where  $h(x)=\frac{\partial}{\partial \theta}g_{\theta}(x)|_{\theta=0}$ and $u(x)=\frac{\partial^2}{\partial \theta^2}g_{\theta}(x)|_{\theta=0}$), and compute the local approximate Bahadur relative efficiencies (LABEs) of the GE test with respect to the LR test. These efficiences are reported in Table \ref{fig: energyBE}.

\begin{table}
\caption{LABE of the energy test for normality}
\centering
\begin{tabular} {c|cccccc}
    $\gamma$ & $g^{(1)}$& $g^{(2)}$ &$g^{(3)}$&$g^{(4)}(1,1)$&$g^{(4)}(0.5,1)$&$g^{(4)}(0,0.5)$  \\\hline
    0.1 &0.501    &0.714&0.843&0.323&0.431&0.630\\
    0.2 &0.520   &0.734&0.861&0.336&0.447&0.636 \\
    0.3 &0.538  &0.754&0.877&0.349&0.464&0.640\\
    0.4 & 0.556   &0.772&0.892&0.362&0.480&0.643\\
    0.5 & 0.573   &0.790&0.906&0.374&0.496&0.645\\
    0.6 & 0.590  &0.806&0.918&0.387&0.512&0.645\\
    0.7 & 0.608   &0.821&0.930&0.399&0.527&0.645\\
    0.8 & 0.623   &0.837&0.940&0.411&0.542&0.643\\
    0.9 & 0.639   &0.851&0.950&0.423&0.558&0.640\\
    1.0 & 0.655&0.865&0.958&0.434&0.572&0.636\\
    1.1 & 0.670   &0.877&0.966&0.446&0.586&0.632\\
    1.2 &0.685    &0.889&0.973&0.457&0.600&0.626\\
    1.3 & 0.699   &0.900&0.978&0.468&0.614&0.620\\
    1.4 & 0.713  &0.911&0.984&0.479&0.628&0.613\\
    1.5 & 0.727   &0.921&0.988&0.490&0.641&0.605\\
    1.6 &0.740  &0.930&0.991&0.501&0.654&0.596\\
    1.7 &0.753&0.939&0.994&0.511&0.666&0.587\\
    1.8 &0.765   &0.947&0.997&0.521&0.679&0.577\\
    1.9 & 0.777   &0.954&0.998&0.531&0.691&0.567\\ \hline
\end{tabular}
\label{fig: energyBE}
\end{table}

\subsubsection{The BHEP test} 
The BHEP test may be  formulated as in \eqref{energy} 
%\begin{equation} \label{bhep}
%T_{n,w}=n\: \int_{-\infty}^\infty \left|\phi_n(t)-e^{-\frac{t^2}{2}}\right|^2 w(t) {\rm{d}}t,
%\end{equation}
with weight function $w(t)=e^{-\gamma t^2}, \ \gamma>0$, which leads to the convenient expression 
\be \label{bhep1}
T_{n,w}=\sqrt{\frac{\pi}{\gamma}} \frac{1}{n} \sum_{j,k=1}^n e^{-\frac{(Z_j-Z_k)^2}{4\gamma}}+n \sqrt{\frac{\pi}{1+\gamma}} -2 \sqrt{\frac{2\pi}{1+2\gamma}} \sum_{j=1}^n e^{-\frac{Z_j^2}{2+4\gamma}}.
\ee%\footnote {This formulation is slightly simpler, but nevertheless equivalent in all essential aspects to the weight function employed by 
%\cite{epps}, and  \cite{henze1997new}, or \cite{tenreiro2009choice}} 
The asymptotic null distribution of the BHEP test  along with the expression 
for the covariance kernel $K(s,t)$ may be found in Henze and Wagner (1997), and corresponding eigenvalues have been recently computed by \cite{ebner}.

%Ebner and Henze (2021b). 

The probability limit of the BHEP test is obtained analogously to \eqref{bTenergy} as
\begin{align*}
 b(\theta)&=\sqrt{\frac{\pi}{1+\gamma}}+
 \sqrt{\frac{\pi}{\gamma}}\mathbb E_{\theta}\left[\exp\Big(-\frac{(X_1-X_2)^2}{4\gamma\sigma^2(\theta)}\Big)\right]
 \\&-2\sqrt{\frac{2\pi}{1+\gamma}}\mathbb E_{\theta}\left[\exp\Big(-\frac{(X_1-\mu(\theta))^2}{(2+4\gamma)\sigma^2(\theta)}\Big)\right],
\end{align*}
while the calculation of efficiencies is also carried out in the same way as in the previous subsection by means of the approximation outlined in \S\ref{sect:eigen}.
These efficiency results (LABEs) for the BHEP test are reported in Table \ref{fig: BHEP}.

\begin{table}
\caption{LABE of the BHEP test for normality}
\centering
\begin{tabular}{c|cccccc}
      $\gamma$ & $g^{(1)}$& $g^{(2)}$ &$g^{(3)}$&$g^{(4)}(1,1)$&$g^{(4)}(0.5,1)$&$g^{(4)}(0,0.5)$  \\\hline
    0.1 &0.477&0.701&0.840&0.302&0.406&0.654\\
    0.2 &0.582&0.814&0.929&0.376&0.501&0.676\\
    0.3 &0.655&0.879&0.968&0.429&0.568&0.658\\
    0.4 &0.710&0.921&0.986&0.471&0.620&0.628\\
    0.5 &0.752&0.948&0.992&0.505&0.661&0.593\\
    0.6 &0.785&0.967&0.993&0.532&0.695&0.559\\
    0.7 &0.812&0.979&0.990&0.555&0.722&0.527\\
    0.8 &0.834&0.987&0.986&0.574&0.745&0.497\\
    0.9 &0.853&0.993&0.980&0.591&0.764&0.469\\
    1 &0.868&0.997&0.974&0.605&0.780&0.443\\
    2 &0.941&0.992&0.917&0.681&0.865&0.281\\
    3 &0.963&0.975&0.879&0.711&0.896&0.202\\
    4 &0.972&0.961&0.855&0.726&0.910&0.158
    \\
    5 &0.977&0.951&0.839&0.735&0.918&0.129\\
    6 &0.979&0.942&0.826&0.741&0.923&0.109\\
    7&0.980&0.936&0.817&0.745&0.926&0.094\\
    8 &0.981&0.931&0.810&0.746&0.929&0.083\\
    9 &0.981&0.926&0.804&0.750&0.930&0.074\\
    10 &0.981&0.923&0.799&0.751&0.932&0.067\\ \hline
    \end{tabular}
\label{fig: BHEP}
\end{table}

\subsubsection{Discussion}

From Tables \ref{fig: energyBE} and \ref{fig: BHEP} we can see that there is a significant influence of the 
tuning parameter on the efficiency of both tests. Specifically, in the case of the GE test, the efficiencies generally grow with $\gamma$,
hence a high value of the tuning parameter (close to the boundary value $\gamma=2$) can be recommended  for this test. On the other hand, the corresponding LABEs of the BHEP test exhibit no such consistent pattern, with the impact on the tuning parameter depending very much on the specific alternative. Specifically for Lehmann and location/scale contamination alternatives higher values of the tuning parameter should be used, while for all other alternatives an ``in--between" value in the neighborhood of $\gamma=1$ yields  better efficiency for the corresponding BHEP test. Moreover the aforementioned value seems to be a good compromise, and if one must choose a single test between the GE and BHEP tests, then the latter test with $\gamma=1$ appears to yield a good overall efficiency.  We can also compare the efficiency of the ECF--based tests to the efficiencies of corresponding tests based on the empirical distribution function (EDF) provided by  \cite{milovsevic2021bahadur}. In this connection, a close inspection of the corresponding efficiency figures shows a significant advantage of the ECF--based tests over their EDF counterparts.

%\subsection{Tests for the Laplace distribution} 
%Recent Monte Carlo studies (see Batsidis et al., 2021, and Desgagn\'e et al., 2021) have shown that the tests suggested by Meintanis (2005) are amongst the most powerful for testing GOF to the Laplace distribution. These tests are formulated as  
%\begin{equation} \label{laplace}
%T_{n,w}=n\: \int_{-\infty}^\infty \left|\phi_n(t)-\frac{1}{1+t^2}\right|^2 w(t) {\rm{d}}t,
%\end{equation}
%with the ECF $\phi_n(t)$ obtained as in \eqref{ECF} with $X_j$ replaced by $x_j=(X_j-\widehat \delta_n)/\widehat c_n$, where $(\widehat \delta_n,\widehat c_n)$ denote the moment or maximum likelihood estimators of the parameters $(\delta,c)$ of the Laplace distribution. The kernel $K(s,t)$ corresponding to each estimation method, as well as explicit expressions for the test statistics with weight function $w(t)=(1+t^2)^2 e^{-\gamma |t|^\beta}, \gamma>0$, and $\beta=1,2$, are provided  by Meintanis (2005). 

\subsection{Tests for the Logistic distribution} 
In complete analogy to the GE and BHEP tests formulated as in \eqref{ts} and \eqref{energy}, \cite{meintanis2004goodness} defines a GOF test statistic for the logistic distribution as   
\begin{equation} \label{logistic}
T_{n,w}=n\: \int_{-\infty}^\infty \left|\phi_n(t)-\frac{\pi t}{\sinh(\pi t)}\right|^2 w(t) {\rm{d}}t,
\end{equation}
with the ECF $\phi_n(t)$ obtained as in \eqref{ECF} with $X_j$ replaced by $Z_j=(X_j-\widehat \mu_n)/\widehat \sigma_n$,
where $(\widehat \mu_n,\widehat \sigma_n)$ denote the moment estimators or the ML estimators of the parameters 
$(\mu,\sigma)$ of the logistic distribution. An explicit test statistic formula corresponding to the weight function  
$w(t)=e^{-\gamma |t|}, \gamma>0$, is given by 
\begin{eqnarray} \label{log1} 
T_{n,w}&=&\frac{2\gamma}{n} \sum_{j,k=1}^n \frac{1}{\gamma^2 +(Z_j-Z_k)^2}-\frac{2}{\pi} \sum_{j=1}^n \left[S^{(1)}_\gamma(Z_j)-\frac{Z^2_j}{2\pi^2} S^{(2)}_\gamma(Z_j)\right] \\ \nonumber &+&\frac{n}{\pi} \left [2 \zeta_2(1+\frac{\gamma}{2\pi})-\frac{\gamma}{\pi}\zeta_3(1+\frac{\gamma}{2\pi}) \right],
\end{eqnarray}
where $\zeta_\gamma(x)=\sum_{k=0}^\infty (k+x)^{-\gamma}$, and
\[
S^{(m)}_\gamma(x)=\sum_{k=0}^\infty\left [ \left(\frac{x}{2\pi}\right)^2+\left(\frac{\gamma+\pi}{2\pi}+k\right)^2 \right]^{-m}.
\]

The expression for the covariance kernel $K(s,t)$ may also be found in \cite{meintanis2004goodness}, along with a Monte Carlo study of the power of the test based on $T_{n,w}$ against the classical tests based on the EDF. These results, nicely complemented by results 
from \cite{gulati2009new}, suggest that the test based on \eqref{log1} is an overall competitive test. %The reader is also referred to \cite{nikitin2020goodness}  for GOF tests based on certain characterizations of the logistic law. 

Turning to efficiency calculations we obtain the corresponding  probability limit (recall \eqref{beta}) as 
\begin{align*}
    b(\theta)&=\mathbb E_\theta\left[\frac{2\gamma}{\gamma^2+\Big(\frac{X_1-X_2}{\sigma(\theta)}\Big)^2}\right]
    -\frac2\pi \mathbb E_\theta \left[S_\gamma^{(1)}\Big(\frac{X-\mu(\theta)}{\sigma(\theta)}\Big)\right]
    \\&-\frac1{\pi^3}\mathbb E_\theta\left[\Big(\frac{X-\mu(\theta)}{\sigma(\theta)}\Big)^2S_\gamma^{(2)}\Big(\frac{X-\mu(\theta)}{\sigma(\theta)}
    \Big)\right]
    +\frac1\pi\Big(2\zeta_2(1+\frac \gamma{2\pi})-\frac \gamma\pi \zeta_3(1+\frac \gamma{2\pi})\Big),
\end{align*}
where $\mu(\theta)$ and $\sigma(\theta)$ denote the probability limits of $\widehat \mu_n$ and $\widehat \sigma_n$, respectively.

In order to calculate $b''(0)$, we use as before numerical integration in
Wolfram Mathematica facilitated by the formulae for the first and second derivatives at zero of $\mu(\theta)$ 
and $\sigma(\theta)$. In this connection and since for the ML estimators there exist no closed expressions for $\mu(\theta)$ and $\sigma(\theta)$, the necessary derivatives  are obtained via the implicit function theorem as  %It can be shown that $b(\theta)=\frac{1}{2}b''(0)\theta^2+o(\theta^2).$
%
%1st case: parameters $\mu$ and $\sigma$ are estimated using MLE
%
%In the case when $\widehat \delta_n$ and $\widehat c_n$ are MLE's, there are no closed 
%expressions for $\mu(\theta)$ and $\sigma(\theta)$, hence the necessary derivatives  are obtained 
%via the implicit function theorem and given in \eqref{miprimlog}, \eqref{sigmaprimlog} and 
%\eqref{midrugolog}-\eqref{sigmadrugolog} given below. 
\begin{align*}  \mu'(0)&=6 \int_{-\infty}^{\infty}\frac{h(x)}{1+e^{-x}}{\rm{d}}x; \; \;
\sigma'(0)=\frac{9}{\pi^2+3}\int_{-\infty}^{\infty}\frac{(1-e^{-x})xh(x)}{1+e^{-x}} {\rm{d}}x,\\
    \mu''(0)&=6\Big(\int \frac{u(x)}{1+e^{-x}}dx + \frac16\mu'(0)\sigma'(0)-2\int \frac{(\mu'(0)+x\sigma'(0))e^{-x}h(x)}{(1+e^{-x})^2}dx \Big),\\
    \sigma''(0)&=\frac{9}{3+\pi^2}\Big(\int xu(x)dx+\frac{1}{2}\mu'(0)^2-4\mu'(0)\sigma'(0)+\frac{1}{18}(\pi^2-6)\sigma'(0)^2
    \\&-4\int\frac{(e^x (x-1)-1)\mu'(0)+e^xx^2\sigma'(0)}{(1+e^x)^2}h(x)dx-2\int \frac{xe^{-x}u(x)}{(1+e^{-x})}\Big).
\end{align*}

On the other hand, when $\widehat \mu_n$ and $\widehat \sigma_n$ are obtained by the method of moments, and by using the corresponding probability limits 
%\begin{align*}
%    \tilde{\mu}& =\bar{X}_n \text{ and}\\
%    \tilde{\sigma}&=\frac{\sqrt{3}}{\pi}\bar{S}_n.
%\end{align*}
\begin{align*}
\mu(\theta)&=\int_{-\infty}^{\infty}xg_\theta(x)dx,\\
\sigma(\theta)&=\frac{\sqrt{3}}{\pi}\bigg(\int_{-\infty}^{\infty}x^2g_\theta(x)dx
-\mu^2(\theta)\bigg)^{\frac{1}{2}},
\end{align*}
we obtain the necessary derivatives as 
\begin{align*}
\mu'(0)&=\int_{-\infty}^{\infty}xh(x)dx; \ \ \mu''(0)=\int_{-\infty}^{\infty}xu(x)dx;\\
%\sigma_{m}'(\theta)&=\frac{\sqrt{3}}{\pi}\frac{\Big(\int_{-\infty}^{\infty}x^2g'(x;\theta)dx-2\int_{-\infty}^{\infty}xg(x;\theta)dx\cdot \int_{-\infty}^{\infty}xg'(x;\theta)dx\Big)}{2\Big(\int_{-\infty}^{\infty}x^2g(x;\theta)dx-(\int_{-\infty}^{\infty}xg(x;\theta)dx)^2\Big)^{\frac{1}{2}}}\\
\sigma'(0)&%=\frac{\sqrt{3}}{\pi}\frac{\Big(\int_{-\infty}^{\infty}x^2h(x)dx\Big)}{2\Big(\frac{\pi^2}{3}\Big)^{\frac{1}{2}}}
=\frac{3}{2\pi^2}\int_{-\infty}^{\infty}x^2h(x)dx;\\
\sigma''(0)%&=\frac{\sqrt{3}}{2\pi}\frac{\Big(\int_{-\infty}^{\infty}x^2u(x)dx-2(\int_{-\infty}^{\infty} xh(x)dx)^2)\cdot\frac{\pi}{\sqrt{3}}-\int_{-\infty}^{\infty}x^2 h(x)dx\cdot\frac{\int_{-\infty}^{\infty}x^2h(x)dx}{2\frac{\sqrt{3}}{\pi}} }{\frac{\pi^2}{3}}
%\\&=\frac{3\sqrt{3}}{2\pi^3}\Big(\frac{\pi}{\sqrt{3}}\int_{-\infty}^{\infty}x^2u(x)dx-(\int_{-\infty}^{\infty}xh(x)dx)^2\frac{2\pi}{\sqrt{3}}-\int_{-\infty}^{\infty}x^2h(x)dx\frac{\pi}{2\sqrt{3}}\Big)
&=\frac{3}{2\pi^2}\int_{-\infty}^{\infty}x^2u(x)dx-
\frac{3}{\pi^2}\Big(\int_{-\infty}^{\infty}xh(x)dx\Big)^2-\frac{3}{4\pi^2}\int_{-\infty}^{\infty}x^2h(x)dx,
\end{align*}
(recall  $h(x)=\frac{\partial}{\partial \theta}g_{\theta}(x)|_{\theta=0}$ and $u(x)=\frac{\partial^2}{\partial \theta^2}g_{\theta}(x)|_{\theta=0}$), and thereby calculate the efficiencies of the test for the logistic distribution based on \eqref{log1}.  These efficiencies (LABEs) are reported in Table \ref{fig: logisticka1}  (ML estimation) and  Table \ref{fig: logisticka2} (moment estimation).

\begin{table} %\label{fig: logt1} 
\caption{LABE of the test  for the logistic distribution (ML estimation)}
\centering
\begin{tabular}{c|cccccc}
      $\gamma$ & $g^{(1)}$& $g^{(2)}$ &$g^{(3)}$&$g^{(4)}(1,1)$&$g^{(4)}(0.5,1)$&$g^{(4)}(0,0.5)$  \\\hline 
     0.1 & 0.274 & 0.456 & 0.641 & 0.468 & 0.463 & 0.702 \\
     0.2 & 0.314 & 0.508 & 0.710 & 0.525 & 0.516 & 0.759 \\
     0.3 & 0.347 & 0.548 & 0.762 & 0.570 & 0.558 & 0.795 \\
     0.4 & 0.378 & 0.581 & 0.804 & 0.607 & 0.593 & 0.820 \\
     0.5 & 0.406 & 0.608 & 0.839 & 0.641 & 0.622 & 0.837 \\
    0.6 & 0.432 & 0.632 & 0.868 & 0.670 & 0.648 & 0.847 \\
    0.7 & 0.457 & 0.652 & 0.893 & 0.695 & 0.669 & 0.852 \\
    0.8 & 0.481 & 0.668 & 0.913 & 0.718 & 0.688 & 0.854 \\
    0.9 & 0.503 & 0.683 & 0.931 & 0.738 & 0.704 & 0.853 \\
    1 & 0.524 & 0.695 & 0.945 & 0.756 & 0.718 & 0.850 \\
    2 & 0.693 & 0.728 & 0.989 & 0.844 & 0.768 & 0.745 \\
 %    2.5 && &&0.848&0.756&0.675\\
 %      2.7 && &&0.845&0.746&0.647\\
 %     2.9 && &&0.840&0.735&0.619\\  
    3 & 0.802 & 0.673 & 0.945 & 0.837 & 0.728 & 0.605 \\
    4 & 0.870 & 0.581 & 0.876 & 0.787 & 0.650 & 0.473 \\
    5 & 0.908 & 0.479 & 0.812 & 0.722 & 0.559 & 0.362 \\
    6 & 0.928 & 0.381 & 0.767 & 0.663 & 0.473 & 0.274 \\
    7 & 0.938 & 0.299 & 0.745 & 0.616 & 0.402 & 0.210 \\
    8 & 0.943 & 0.237 & 0.739 & 0.584 & 0.351 & 0.166 \\
    9 & 0.946 & 0.194 & 0.741 & 0.565 & 0.315 & 0.136 \\
    10 & 0.948 & 0.163 & 0.744 & 0.552 & 0.230 & 0.116 \\ \hline
    \end{tabular}
\label{fig: logisticka1}
\end{table}

\begin{table} 
\caption{LABE of the test  for the logistic distribution (moment estimation)}
\centering
\begin{tabular}{c|cccccc}
     $\gamma$ & $g^{(1)}$& $g^{(2)}$ &$g^{(3)}$&$g^{(4)}(1,1)$&$g^{(4)}(0.5,1)$&$g^{(4)}(0,0.5)$  \\\hline
     0.1 & 0.485 & 0.695 & 0.811 & 0.680 & 0.693 & 0.878\\
     0.2 & 0.539 & 0.754 & 0.874 & 0.741 & 0.752 & 0.924\\
     0.3 & 0.582 & 0.794 & 0.915 & 0.786 & 0.794 & 0.947\\
     0.4 & 0.619 & 0.824 & 0.944 & 0.820 & 0.825 & 0.957\\
     0.5 & 0.651 & 0.846 & 0.965 & 0.846 & 0.849 & 0.958\\
    0.6 & 0.670 & 0.863 & 0.980 & 0.868 & 0.867 & 0.954 \\
    0.7 & 0.706 & 0.875 & 0.989 & 0.884 & 0.880 & 0.945\\
    0.8 & 0.729 & 0.883 &  0.995 & 0.890 & 0.890 & 0.933\\
    0.9 & 0.751 & 0.889 & 0.998 & 0.908 & 0.896 & 0.919\\
    1 & 0.770 & 0.891 & 0.998 & 0.916 & 0.901 & 0.904\\
    2 & 0.898 & 0.841 & 0.931 & 0.910 & 0.863 & 0.724\\
 %    2.5 && &&0.848&0.756&0.675\\
 %      2.7 && &&0.845&0.746&0.647\\
 %     2.9 && &&0.840&0.735&0.619\\  
    3 & 0.941 & 0.729 & 0.816 & 0.828 & 0.760 & 0.559 \\
    4 & 0.773 & 0.504 & 0.580 & 0.601 & 0.532 & 0.355 \\
    5 & 0.581 & 0.323 & 0.387 & 0.405 & 0.348 & 0.213 \\
    6 & 0.422 & 0.202 & 0.254 & 0.266 & 0.222 & 0.126\\
    7 & 0.304 & 0.127 & 0.169 & 0.176 & 0.142 & 0.076\\
    8 & 0.221 & 0.081 & 0.114 & 0.118 & 0.092 & 0.047\\
    9 & 0.161 & 0.052 & 0.079 & 0.080 & 0.061 & 0.029\\
    10 & 0.120 & 0.034 & 0.055 & 0.056 & 0.041 & 0.019\\ \hline
    \end{tabular}
\label{fig: logisticka2}
\end{table}

For comparison purposes we also provide in Table \ref{tab:logisticEDF} corresponding  efficiencies for three EDF-based tests, namely the Kolmogorov-Smirnov (KS), the Cram\'er-von Mises (CM), and the Anderson-Darling (AD) test. The calculations were carried out using the same method used by  \cite{milovsevic2021bahadur} for normality tests, while the covariance function of the corresponding empirical process is available from \cite{stephens1979tests}.

\begin{table}
    \centering
    \begin{tabular}{c|cccccc}
    Test  & $g^{(1)}$& $g^{(2)}$ &$g^{(3)}$&$g^{(4)}(1,1)$&$g^{(4)}(0.5,1)$&$g^{(4)}(0,0.5)$  \\\hline
        KS &0.191&0.312 &0.455&0.329&0.318& 0.435 \\
        CM &0.475&0.677&0.920&0.715&0.693&0.866\\
        AD&0.754&0.640& 0.999&0.808&0.697&0.718 \\ \hline
    \end{tabular}
    \caption{LABE of EDF tests for the logistic distribution}
    \label{tab:logisticEDF}
\end{table}

\subsubsection{Discussion}
{
We note that this is the first time that efficiency results for a test for the composite hypothesis of GOF to the logistic distribution are obtained. In particular, there exist no results for the Bahadur efficiencies for classical GOF tests in the case of the logistic distribution, or
in fact for any GOF test for this distribution in the case of the ``fully" composite null hypothesis, i.e. when both the location and scale parameters are unknown.} % This is the main reason for inclusion of EDF based tests in this study. }

When ECF based tests are compared among themselves,  from Tables \ref{fig: logisticka1} and \ref{fig: logisticka2} it may be inferred that the efficiencies vary considerably between estimation methods, and that the variant of the test which uses the method  of moments is somewhat more efficient than its counterpart based on ML estimation, a finding that is in line with the Monte Carlo results of Meintanis (2004). This phenomenon of (generally) less efficient estimators, such as the moment estimators, producing more powerful tests has already been pointed out by \cite{Gurtler}  under similar circumstances, and dates back at least to \cite{Kallenberg}. Turning to the tuning  parameter we observe that its impact on test performance is noticeable, and that, by taking into account all considered alternatives, we can recommend $\gamma=3$ in the case of ML estimation and $\gamma=1$ in the case of moment estimation as overall good choices.  

{
As far as classical tests are concern, we notice from Table \ref{tab:logisticEDF} that the AD test is more efficient than the other EDF--based tests and competitive with the ECF--based tests. Nevertheless the moment--based ECF test appears to be preferable overall against the AD in the neighborhood of the suggested value $\gamma=1$. The other EDF--based tests are also less efficient than the ECF tests for the recommended value of the tuning parameter $\gamma$.
}

%Kad se porede ova dva, posotje razlike i metod ocenjivanja utice na efiksanost. MoM je nesto bolji od MLE.
%Nema rezultata za klsicne testove u literaturi za slozenu hioptezu pa ne mozemo da uporedimo. Sto se tice $gamma$,
%o njega zavisi, MoM preporucujemo 1,  MLE 3.

\subsection{Tests for the exponential distribution} 
In extensive Monte Carlo studies (see for instance \cite{henze2005recent} and \cite{grane3}), certain tests based on the ECF
are found to be competitive for testing GOF to the exponential distribution. 
One such test is also a weighted L2--type test, but the corresponding test statistic admits a slightly 
different formulation from that in \eqref{ts}. Specifically we have   
\begin{equation} \label{exptest}
T_{n,w}=n\: \int_{-\infty}^\infty \left(|\phi_n(t)|^2-C_n(t)\right)^2 w(t) {\rm{d}}t,
\end{equation}
with $|\phi_n(t)|$ denoting the modulus and $C_n(t)$ the real part of  the ECF obtained from \eqref{ECF} 
by replacing $X_j$  by 
\be \label{stand1}  Z_j=X_j/\bar X_n, \ j=1,...,n.\ee
As weight functions \cite{henze2005recent} suggest $w_\beta(t)=e^{-\gamma |t|^\beta}, \gamma>0$, $\beta=1,2$, 
and provide explicit expressions  for the resulting test statistics. The corresponding kernel $K(s,t)$ figuring in 
\eqref{eigen} may be computed as 
\[
K(s,t)=\frac{s^2t^2(1+s^2+t^2)}{(1+s^2)(1+t^2)(1+(s-t)^2)(1+(s+t)^2)}.
\]

For the weight function $w_1(t)=e^{-\gamma|t|}$, the probability limit figuring in \eqref{beta}  is given by  
\begin{align*}
 b_1(\theta)&=\gamma \mathbb E_{\theta}\left[\frac{1}{\gamma^2+(\frac{X_1-X_2}{\sigma(\theta)})^2}+ 
 \frac{1}{\gamma^2+(\frac{X_1+X_2}{\sigma(\theta)})^2}\right]\\&-
 2\gamma \mathbb E_{\theta}\left[\frac{1}{\gamma^2+(\frac{X_1-X_2-X_3}{\sigma(\theta)})^2}+ 
 \frac{1}{\gamma^2+(\frac{X_1-X_2+X_3}{\sigma(\theta)})^2}\right]
 \\&+\gamma \mathbb E_\theta\left[\frac{1}{\gamma^2+(\frac{X_1-X_2-X_3+X_4}{\sigma(\theta)})^2}+ 
 \frac{1}{\gamma^2+(\frac{X_1-X_2+X_3-X_4}{\sigma(\theta)})^2}\right],
\end{align*}
while for $w_2(t)=e^{-\gamma t^2}$, the same probability limit is given by 
 \begin{align*}
 b_2(\theta)&=\frac12\sqrt\frac\pi\gamma \: \mathbb E_{\theta}\left[\exp\Big(\frac{(X_1-X_2)^2}{4\gamma\sigma^2(\theta)}\Big)+ 
 \Big(\exp\Big(\frac{(X_1+X_2)^2}{4\gamma\sigma^2(\theta)}\Big)\right]\\&-
 \sqrt\frac\pi\gamma \: \mathbb E_{\theta}\left[\exp\Big(\frac{(X_1-X_2-X_3)^2}{4\gamma\sigma^2(\theta)}\Big)+ 
 \Big(\exp\Big(\frac{(X_1-X_2+X_3)^2}{4\gamma\sigma^2(\theta)}\Big)\right]
\\&+\frac12\sqrt\frac\pi\gamma \: \mathbb E_{\theta}\left[\exp\Big(\frac{(X_1-X_2-X_3+X_4)^2}{4\gamma\sigma^2(\theta)}\Big)+ 
 \Big(\exp\Big(\frac{(X_1-X_2+X_3-X_4)^2}{4\gamma\sigma^2(\theta)}\Big)\right],
\end{align*}
where, recall, $\sigma(\theta)=\int_{-\infty}^\infty xg_\theta(x)dx$.

As far as alternatives are concerned, the exponential distribution has some common close alternatives. Specifically we consider the following alternatives:
\begin{itemize}
		\item the Weibull distribution with  density
		\begin{equation*}
		g_\theta(x)=e^{-x^{1+\theta}}(1+\theta)x^\theta,\theta>0,x\geq0;
		\end{equation*}
		\item the gamma distribution with density
		\begin{equation*}
		g_\theta(x)=\frac{x^\theta e^{-x}}{\Gamma(\theta+1)},\theta>0,x\geq0;
		\end{equation*}
		\item 
		the Makeham distribution with density 
		\begin{equation*}
		g_\theta(x)=e^{-x-\theta e^{x}}(1+\theta e^x),\theta>0,x\geq0;
		\end{equation*}
		\item the linear failure rate (LFR) distribution with density 
		\begin{equation*}
		g_\theta(x)=e^{-x-\theta\frac{x^2}{2}}(1+\theta x),\theta>0,x\geq0;
		\end{equation*}
		\item the mixture of exponential distributions with negative weights (ME($\beta$)) with density
		\begin{equation*}
		g_\theta(x;\beta)=(1+\theta)e^{-x}-\theta\beta e^{-\beta x},\theta\in\left(0,\frac{1}{\beta-1}\right],x\geq0;
		\end{equation*}
	\end{itemize}
These alternatives are also used by e.g. \cite{Publ,nikitin2016efficiency, milovsevic2016asymptotic,cuparic2018new}.  The resulting efficiencies are calculated by following analogous steps as in the previous sections and  are reported  in Table \ref{fig: exp1} and Table \ref{fig: exp2}, for the test based on $w_1(t)=e^{-\gamma |t|}$ and $w_2(t)=e^{-\gamma |t|^2}$, respectively.
  
\begin{table}
\centering
\caption{LABE of the test for exponentiality with weight function $w_1$}
\begin{tabular}{c|cccccc|}
    $\gamma$ & Weibull&Gamma&Makeham&LFR&ME(3)&ME(6)  \\\hline
    0.1&0.502&0.619&0.404&0.125&0.585&0.831\\
    0.2&0.552&0.625&0.470&0.155&0.669&0.841\\
    0.3&0.584&0.625&0.523&0.183&0.724&0.834\\
    0.4&0.610&0.632&0.580&0.203&0.763&0.819\\
    0.5&0.625&0.627&0.621&0.231&0.792&0.802\\
    0.6&0.637&0.626&0.657&0.252&0.813&0.784\\
    0.7&0.649&0.626&0.688&0.275&0.829&0.767\\
    0.8&0.659&0.624&0.713&0.295&0.842&0.750\\
    0.9&0.671&0.622&0.737&0.315&0.852&0.735\\
    1 &0.676 &0.620&0.757&0.334&0.859&0.720\\
    2&0.722&0.596&0.881&0.494&0.874&0.616\\
    3&0.738&0.573&0.925&0.611&0.850&0.553\\
    4&0.741&0.553&0.938&0.699&0.821&0.510\\
    5&0.737&0.534&0.936&0.765&0.792&0.478\\
    6&0.730&0.518&0.928&0.815&0.766&0.453\\
    7&0.721&0.504&0.917&0.853&0.743&0.434\\
    8&0.713&0.487&0.905&0.882&0.723&0.418\\
    9&0.704&0.478&0.893&0.905&0.705&0.404\\
    10&0.696&0.461&0.881&0.923&0.690&0.393\\ \hline
\end{tabular}
\label{fig: exp1}
\end{table}

\begin{table}
\centering
\caption{LABE of the test for exponentiality with weight function $w_2$}
\begin{tabular}{c|cccccc|}
    $\gamma$ & Weibull&Gamma&Makeham&LFR&ME(3)&ME(6)  \\\hline
    0.1&0.581& 0.585 & 0.580 & 0.198 & 0.756 & 0.749\\
    0.2&0.619 & 0.584 & 0.676 & 0.261 & 0.807 & 0.698\\
    0.3&0.639 & 0.581 & 0.735 & 0.312 & 0.828 & 0.665\\
    0.4&0.657 & 0.579 & 0.773 & 0.350 & 0.838 & 0.640\\
    0.5&0.670 &0.576 & 0.802 & 0.384 & 0.843 & 0.621\\
    0.6&0.679 & 0.573 & 0.824 & 0.413 & 0.845 & 0.605\\
    0.7&0.687 & 0.571 & 0.842 & 0.439 & 0.845 & 0.592\\
    0.8&0.693 & 0.568 & 0.856 &0.463 & 0.844 & 0.580\\
    0.9&0.698 & 0.566 & 0.867 & 0.485 & 0.842 & 0.570\\
    1 &0.703 & 0.564 & 0.877 & 0.504 & 0.840 & 0.562\\
    2&0.722 & 0.544 & 0.918 & 0.643 & 0.811 & 0.505\\
    3&0.725 & 0.529 & 0.925 & 0.724 & 0.785 & 0.474\\
    4&0.723 & 0.516 & 0.921 & 0.778 & 0.763 & 0.452\\
    5&0.719 & 0.505 & 0.915 & 0.816 & 0.746 & 0.436\\
    6&0.714 &0.497 & 0.908 & 0.844 & 0.731 & 0.424\\
    7&0.710 &0.489 &0.900 &0.866 & 0.718 & 0.414\\
    8&0.705 & 0.481& 0.892 & 0.892 & 0.708 & 0.406\\
    9&0.700 & 9.476 & 0.885 & 0.885 & 0.698 & 0.399\\
    10&0.697&0.467 & 0.879 & 0.910 & 0.690 & 0.393 \\ \hline
\end{tabular}
\label{fig: exp2}
\end{table}

\subsubsection{Discussion}

From Tables \ref{fig: exp1} and \ref{fig: exp2} one can notice that the type of monotonicity  of the efficiencies with respect to the tuning parameter $\gamma$ varies amongst alternatives. However, for any fixed alternative, both tests exhibit the same behaviour, and if we compare maximal reached efficiencies, the test
with weight function $w_1(\cdot)$ appears to have a slight edge. Next we compare the efficiencies obtained here with the efficiencies reported by \cite{revista} for a wide variety of exponentiality tests and for four of the alternatives considered herein (e.g., Weibull, Gamma, LFR and ME(3)). Specifically the tests based on $T_{n,w}$ seem to be less efficient than some of the best tests considered by \cite{revista}. When restricting comparison to the classical tests based on the EDF, we see that the new tests are more efficient than the KS test, but overall less efficient than the  CM and AD tests, with the exception of the LFR alternative where, with proper choice of $\gamma$, our tests compete well with these two tests.

\section{Tests for the bivariate normal distribution} \label{sec5}

We consider the problem of testing the null hypothesis ${\cal{H}}_0$ that the sample comes from a bivariate normal distribution, using the BHEP and GE tests presented in Section 4.1. Efficiencies are calculated for two cases: The simple hypothesis of known mean/covariance, that the law of $X$ is ${\cal{N}}_2(0,{\rm{I}})$, and the hybrid case of unknown mean/known covariance, i.e. that the law of $X$ is ${\cal{N}}_2(\mu,{\rm{I}})$, where ${\cal{N}}_p(\mu, \Sigma)$ denotes  the multivariate normal distribution with mean vector $\mu$ and covariance matrix equal to $\Sigma$, and ${\rm{I}}$ stands for the identity matrix in the corresponding dimension. The covariance kernels of the test statistics  are available from \cite{mori2021energy} (GE test) and \cite{henze1997new} (BHEP test), while for the simple hypothesis case with $\gamma=1/2$, eigenvalues of the BHEP test in arbitrary dimension may be obtained from \cite{baring}.\footnote{As a numerical confirmation, we have compared the efficiencies obtained by using Baringhaus' eigenvalues with those using the Monte Carlo method and found them to be quite close in all cases} 
%1) We have used the value from the paper you sent us. However, we notice
%that there are some discrepancies in the value we obtained using MC
%simulations. At first, we thought it is due to the large variance of the
%MC estimator, therefore we decided to examine the distribution of
%eigenvalues obtained for fixed N. Our conclusion was that maybe the
%authors made mistake, but we decided to keep this value until we check
%completely their proof. I am highlighting that the differences of obtained
%efficiencies (using their and MC method for the calculation of eigenvalue 
%are not large).  We shall check this once again during the week.

The method of calculation of efficiencies is the same as for univariate tests. Specifically, consider an alternative with density $g(x,y;\theta)$, and for small $\theta$, compute the double KL distance of this alternative to  the null hypothesis,   as
\begin{itemize}
    \item simple hypothesis case: \begin{align*}
        2K(\theta)=\left(\iint_{\mathbb R^2}\frac{h^2(x,y)}{f_0(x,y)}{\rm{d}}x{\rm{d}}y\right) \theta^2+o(\theta^2),
    \end{align*}
    \item estimated mean case:
         \begin{align*}
        2K(\theta)=\left(\iint_{\mathbb R^2}\frac{\big(h(x,y)-f_0(x,y)(x\mu'_X(0)+y\mu'_Y(0))\big)^2}{f_0(x,y)}{\rm{d}}x{\rm{d}}y\right) \theta^2+o(\theta^2),
    \end{align*}
    where \begin{align*}
        \mu_X(\theta)=\iint_{\mathbb R^2} xg(x,y;\theta){\rm{d}}x{\rm{d}}y; \ \ 
        \mu_Y(\theta)=\iint_{\mathbb R^2}yg(x,y;\theta){\rm{d}}x{\rm{d}}y,
    \end{align*}
\end{itemize}
and $f_0(\cdot,\cdot)$ denotes the density of the standard normal distribution. 
The proof is analogous to that of Theorem \ref{TKL}, so we omit it here.

The multidimensional version of the eigenvalue approximation procedure described in Section \ref{sec3} is computationally very complex, hence here we obtained the largest eigenvalues using the Monte Carlo procedure described in \cite{mori2021energy}. In this connection we have noticed that the variance of the Monte Carlo estimator of the eigenvalue is small in the case of BHEP test, which makes the obtained values reliable. On the other hand, for the GE test, the variance is significant, making the approximation not so accurate. In addition, this problem drastically increases with the increase of the tuning parameter $\gamma$, therefore we excluded the cases when $\gamma>1$ from the study.

%Consider the mixture alternative

%\begin{align*}
 %   g(\boldsymbol{x};\theta,\boldsymbol{\mu},\Sigma)=(1-\theta)f(\boldsymbol{x};0,I)+\theta f(\boldsymbol{x};\boldsymbol{\mu},\Sigma), 
%\end{align*}
%where
%\begin{align*}
 %   f(\boldsymbol{x};\boldsymbol{\mu},\Sigma)=\frac{1}{2\pi\sqrt{\det\Sigma}}\exp\Big(-\frac12((\boldsymbol{x}-\boldsymbol{\mu})^T\Sigma^{-1}(\boldsymbol{x}-\boldsymbol{\mu})\Big),
%\end{align*}
%$\boldsymbol{\mu}$ is the vector of means, and $\Sigma$ is the variance-covariance matrix.

%The Kullback-Leibler distance of $g(\boldsymbol{x};\theta,\boldsymbol{\mu},\Sigma)$ to the set of all 
%bivariate normal distributions is

%\begin{align*}
 %   2K(\theta)&=\int \frac{(f(\boldsymbol{x};\boldsymbol{\mu},\Sigma)-f(\boldsymbol{x};0,I))^2}{f(\boldsymbol{x};0,I)} d\boldsymbol{x} - ...\\
%\end{align*}

%general case
%\begin{align*}
 %   2K(\theta)&=\int \frac{h(\boldsymbol{x};\boldsymbol{\mu},\Sigma)^2}{f(\boldsymbol{x};0,I)} d\boldsymbol{x} \cdot\theta^2\\
%\end{align*}
In order to introduce our alternatives write $f(\mu_1,\mu_2,\sigma_1,\sigma_2,\rho)$ for the bivariate normal density with mean vector $\mu=(\mu_1,\mu_2)^\top$ and covariance matrix $\Sigma$ with elements $\sigma_{11}=\sigma^2_1$, $\sigma_{22}=\sigma^2_2$, and $\sigma_{12}=\sigma_{21}=\rho$. With this notation we consider the following alternatives for the case of the simple hypothesis tests:  \begin{itemize}
\item location alternative $g_{l}(\theta)=f(\theta,0,1,1,0)$
\item correlation alternative
$g_{c}(\theta)=f(0,0,1,1,\theta)$
\item single scale alternative
$g_{s1}(\theta)=f(0,0,1-\theta,1,0)$
\item double scale alternative
$g_{s2}(\theta)=f(0,0,1-\theta,1-\theta,0)$
\item contamination alternative
$$g_{cn}(\mu_1,\mu_2,\sigma_1,\sigma_2,\rho,\theta)=(1-\theta)f(0,0,1,1,0)+\theta f(\mu_1,\mu_2,\sigma_1,\sigma_2,\rho),$$
for the specific choices of parameters given in Table \ref{tab: contamination}.
\end{itemize}
\begin{table}
\centering
\caption{Parameters of the contaminating distribution}
\begin{tabular}{c|ccccc}
 &  $\mu_1$&$\mu_2$&$\sigma_1$&$\sigma_2$&$\rho$ \\\hline
  $g^{(1)}_{cn}$ & 0.1 & 0 & 1 & 1 &0\\
  $g^{(2)}_{cn}$ & 0.5 & 0 & 1 & 1 &0\\
   $g^{(3)}_{cn}$ & 0.9 & 0 & 1 & 1 &0\\
    $g^{(4)}_{cn}$ & 1.5 & 0 & 1 & 1 &0\\\hline
$g^{(5)}_{cn}$ & 0 & 0 & 1 & 1 &0.1\\
$g^{(6)}_{cn}$ & 0 & 0 & 1 & 1 &0.5\\
$g^{(7)}_{cn}$ & 0 & 0 & 1 & 1 &0.9\\\hline
$g^{(8)}_{cn}$ & 0 & 0 & 0.5 & 1 &0\\
$g^{(9)}_{cn}$ & 0 & 0 & 0.7 & 1 &0\\
$g^{(10)}_{cn}$ & 0 & 0 & 0.9 & 1 &0\\
$g^{(11)}_{cn}$ & 0 & 0 & 1.1 & 1 &0\\\hline
\end{tabular}
\label{tab: contamination}
\end{table}

\subsection{Discussion}
Table \ref{tab: LABEbivariateSimple} (simple hypothesis case) and Table \ref{tab: LABEunknownM} (estimated mean case)\footnote{For the hybrid case  we only present results for contamination alternatives as location alternative are clearly excluded and  efficiencies for the correlation and scale alternatives coincide with those of the simple hypothesis case} report efficiency results for the GE and BHEP tests for bivariate normality. To facilitate comparison, for each alternative we report in italics the best efficiency of the test under discussion against  $\gamma$, while a bold entry indicates the best efficiency value across both tests.    

For the GE test the monotonicity is clear, since better efficiencies are obtained as $\gamma$ increases, in both the simple and the hybrid case test, a behaviour that mimics the univariate tests. On the other hand, also analogously to the univariate case,  no such pattern is visible for the BHEP test in the simple hypothesis case, as efficiencies are better at the one or the other end of the tuning parameter interval, or even in--between values of $\gamma$, depending on the type of alternative. Things are somewhat more clear though for the BHEP test with estimated mean, and specifically in this case better efficiencies are observed for larger values of $\gamma$, nearly uniformly over alternatives.    

When comparing the GE and BHEP tests with each other, one cannot suggest a single test as being uniformly more efficient over all alternatives. Nevertheless in the simple hypothesis case, the GE test with $\gamma=1$  appears to perform better against location contamination alternatives, while the BHEP test with $\gamma=0.25, 0.50$ shows higher efficiency in nearly all other cases of alternatives. With a few exceptions, the BHEP test with $\gamma=2$ seems to be also preferable over the GE test in the case of testing for bivariate normality with unknown mean.

\begin{table}
\centering
\caption{LABE for the GE and BHEP tests--simple hypothesis}
\begin{tabular}{c|ccc|cccc|}
Test & \multicolumn{3}{|c|}{GE} &\multicolumn{4}{|c|}{BHEP}\\\hline
    $\gamma$ &  0.5 & 0.7 & 1  & 0.25 & 0.5 & 1 & 2\\\hline
$g_{l}$    & 0.285& 0.564 &\bf 0.961 &0.631&0.762 &0.860&\it 0.927 \\
$g_{c}$&0.054 & 0.092 & \it 0.119 &0.252&\bf 0.254&0.215&0.154\\
$g_{s1}$&0.080 & 0.138 &\it 0.180 &0.379&\bf 0.381&0.323&0.232\\
$g_{s2}$&0.106 & 0.184 & \it 0.240&0.505&\bf 0.508&0.430&0.309\\\hline
$g^{(1)}_{cn}$& 0.284 & 0.562& \bf 0.956 &0.627&0.735&0.853&\it 0.918\\
$g^{(2)}_{cn}$& 0.248& 0.492&\bf 0.840 &0.542&0.637&0.743&\it 0.808\\
$g^{(3)}_{cn}$& 0.178& 0.355 &\bf 0.609 &0.378&0.449&0.531&\it 0.582\\
$g^{(4)}_{cn}$& 0.068& 0.137 &\bf 0.238 &0.134&0.163&0.200&\it 0.225\\\hline
$g^{(5)}_{cn}$& 0.053& 0.092 &\it 0.119 &\bf 0.253&0.245&0.213&0.150\\
$g^{(6)}_{cn}$& 0.046& 0.078 &\it 0.101 &\bf 0.237&0.215&0.172&0.117\\
$g^{(7)}_{cn}$& 0.019& 0.031 &\it 0.038 &\bf 0.110&0.082&0.055&0.033\\\hline
$g^{(8)}_{cn}$& 0.098& 0.159 &\it 0.193&\bf 0.588&0.454&0.310&0.180\\
$g^{(9)}_{cn}$& 0.102& 0.170 &\it 0.213 &\bf 0.563&0.486&0.369&0.233\\
$g^{(10)}_{cn}$& 0.090& 0.156 &\it 0.199 &\bf 0.451&0.425&0.355&0.245\\
$g^{(11)}_{cn}$& 0.067& 0.116 &\it 0.154 &\bf 0.303&0.303&0.275&0.206\\ \hline
\end{tabular}   
\label{tab: LABEbivariateSimple}
\end{table}
\bigskip

\begin{table}
\centering
    \caption{LABE for the GE and BHEP tests -- estimated mean}
\begin{tabular}{c|ccc|cccc|}
Test & \multicolumn{3}{|c|}{GE} &\multicolumn{4}{|c|}{BHEP}\\\hline
    $\gamma$ &  0.5 & 0.7 & 1  & 0.25 & 0.5 & 1 & 2\\\hline
$g^{(1)}_{cn}$   & 0.218& 0.366 &\it 0.641 &0.455&0.553 &0.637&\bf 0.700 \\
$g^{(2)}_{cn}$& 0.200& 0.330&\it 0.563 &0.409&0.499&0.577&\bf 0.636\\
$g^{(3)}_{cn}$& 0.155& 0.257 &\it 0.440 &0.313&0.386&0.453&\bf 0.505\\
$g^{(4)}_{cn}$& 0.075& 0.126 &\it 0.217 &0.146&0.184&0.222&\bf 0.258\\\hline
$g^{(5)}_{cn}$& 0.147& 0.243 &\it 0.412 &0.305&0.369&0.423&\bf 0.465\\
$g^{(6)}_{cn}$& 0.128& 0.208 &\it 0.348 &0.285&0.323&0.348&\bf 0.364\\
$g^{(7)}_{cn}$& 0.054& 0.083 &\it 0.130 &\bf 0.137&0.124&0.111&0.103\\\hline
$g^{(8)}_{cn}$& 0.269& 0.422 &\it 0.667 &\bf 0.710&0.683&0.614&0.546\\
$g^{(9)}_{cn}$& 0.281& 0.451 &\bf 0.736 &0.680&0.731&\it 0.731&0.704\\
$g^{(10)}_{cn}$& 0.250& 0.409 &\it 0.687 &0.545&0.639&0.704&\bf 0.740\\
$g^{(11)}_{cn}$& 0.185& 0.309 &\it 0.531 &0.363&0.455&0.544&\bf 0.624\\ \hline
\end{tabular}
\label{tab: LABEunknownM}
\end{table}

\section{Conclusion} \label{sec6} 
{We consider test optimality in the Bahadur sense for certain ECF--based GOF tests for the normal, the logistic and the exponential distribution. In the case of testing for normality we compare the efficiencies of the celebrated generalized energy and BHEP tests, and we offer suggestions as to which tests should be used on the basis of our efficiency comparisons, in the univariate as well as the bivariate case, with or without estimated parameters. In the case of testing for the logistic distribution efficiency comparisons of an ECF--based test against classical tests  based on the EDF are reported for moment as well as maximum likelihood estimation of parameters, whereas for testing exponentiality we compare our efficiency findings with those of \cite{revista} for a wide range of alternative tests, including classical ones. Overall ECF tests appear to  compare well and often outperform competitors. Moreover  as ECF--based tests involve a  tuning parameter, the results reported herein may be used in determining which value of this parameter should be employed by the user, a problem that in various forms occupies researchers to this date; see for instance the contributions of   \cite{ebner2021bahadur}, \cite{ebner}, \cite{Tenreiro19}, and \cite{allison}.}         

%The choice of these families comprise  some of the most popular distributional models, 
%and have also been considered elsewhere in analogous optimality studies; see for instance \cite{dewet}, 
%and references therein.     

%\medskip

%\noindent {\bf Acknowledgement}   
%The work of B. Milo\v sevi\'c is supported by the Ministry of Education, Science and Technological Development of Republic of Serbia.
\section*{Appendix} 
\textbf{Proof of Theorem \ref{TKL}.} Let $\mu=\mu(\theta)$ and $\sigma=\sigma(\theta)$ be the values that minimize \eqref{KL}. 
Then, differentiating \eqref{KL} with respect to $\theta$ we get
%\begin{align*}
%K(\theta)%&=\int\log g(x,\theta) g(x,\theta)dx-\int\log \frac{1}{\sigma(\theta)}f(\frac{x-m(\theta)}{\sigma(\theta)}) \cdot g(x,\theta)dx\\
%&=\int\log g(x,\theta) g(x,\theta)dx +\log \sigma(\theta)-\int\log f(\frac{x-\mu(\theta)}{\sigma(\theta)}) 
%\cdot g(x,\theta)dx
%\end{align*}
\begin{align*}
K'(\theta)&=\int \log g_\theta(x)g_\theta'(x)dx-
\int \frac{1}{f(\frac{x-\mu(\theta)}{\sigma(\theta)})}f'\Big(\frac{x-\mu(\theta)}{\sigma(\theta)}\Big) 
\Big(\frac{x-\mu(\theta)}{\sigma(\theta)}\Big)'g_\theta(x)
\\&+\frac{\sigma'(\theta)}{\sigma(\theta)}-\int\log f\Big(\frac{x-\mu(\theta)}{\sigma(\theta)}\Big) 
g_\theta'(x)dx.
\end{align*}
It is easy to show that $K'(0)=0$. Differentiating \eqref{KL} once more we obtain at $\theta=0$, 
%\begin{align*}
%    K''(\theta)&= \int \frac{(g'(x;\theta))^2}{g(x;\theta)}dx + \int \log g(x;\theta)g''(x;\theta) dx + 
%    \frac{\sigma''(\theta)\sigma(\theta)-(\sigma'(\theta))^2}{\sigma^2(\theta)} 
%    \\&+ \int \frac{1}{f^2(\frac{x-\mu(\theta)}{\sigma(\theta)})}f'(\frac{x-\mu(\theta)}{\sigma(\theta)}) 
%    \big(\frac{x-\mu(\theta)}{\sigma(\theta)})'\big)^2g(x,\theta)dx 
%    \\&- \int \frac{1}{f(\frac{x-\mu(\theta)}{\sigma(\theta)})}f''(\frac{x-\mu(\theta)}{\sigma(\theta)}) 
%    \big((\frac{x-\mu(\theta)}{\sigma(\theta)})'\big)^2g(x,\theta)dx 
%    \\&- \int \frac{1}{f(\frac{x-\mu(\theta)}{\sigma(\theta)})}f'(\frac{x-\mu(\theta)}{\sigma(\theta)}) 
%    (\frac{x-\mu(\theta)}{\sigma(\theta)})''g(x;\theta)dx
%    \\&-\int \frac{1}{f(\frac{x-\mu(\theta)}{\sigma(\theta)})}f'(\frac{x-\mu(\theta)}{\sigma(\theta)}) 
%    (\frac{x-\mu(\theta)}{\sigma(\theta)})'g'(x,\theta)dx 
%    \\&- \int \frac{f'(\frac{x-\mu(\theta)}{\sigma(\theta)})
%    (\frac{x-\mu(\theta)}{\sigma(\theta)})'}{f(\frac{x-\mu(\theta)}{\sigma(\theta)})} g'(x;\theta) dx - 
%    \int \log f(\frac{x-\mu(\theta)}{\sigma(\theta)}) g''(x;\theta)dx 
%\end{align*}
%Putting $\theta=0$, 
\begin{align*}
    K''(0)&=\int\frac{h^2(x)}{f(x)}dx+\int \log f(x) u(x)dx+\sigma''(0)-(\sigma'(0))^2
    \\&+\int\frac{1}{f(x)}(f'(x))^2(\mu'(0)+x\sigma'(0))^2dx-\int f''(x)(\mu'(0)+x\sigma'(0))^2dx
%    \end{align*}
    \\   &-\int f'(x)(-\mu''(0)+2 \mu'(0) \sigma'(0)+x \left(2 \sigma'(0)^2-\sigma''(0)\right))dx \\
%  \end{align*}  
%  \begin{align*}
    &+\int\frac{1}{f(x)}f'(x)(\mu'(0)+x\sigma'(0))h(x)dx\\
    &+\int\frac{1}{f(x)}f'(x)(\mu'(0)+x\sigma'(0))h(x)dx-\int \log f(x)u(x)dx,
\end{align*}
where $h(x)=\frac{\partial}{\partial \theta}g_{\theta}(x)|_{\theta=0}$ and $u(x)=\frac{\partial^2}{\partial \theta^2}g_{\theta}(x)|_{\theta=0}$. Rearranging  terms in the above expression, and expanding $K(\theta)$ into a Maclaurin expansion completes  
the proof.
 \hfill$\Box$

\bibliographystyle{chicago}
\bibliography{biblio}

\end{document}